\documentclass{article}
\usepackage{amsfonts,amsmath}
\mathchardef\mhyphen="2D 
\usepackage[a4paper]{geometry}
\usepackage{microtype}
\usepackage[shortlabels]{enumitem} 
\usepackage{relsize}

\def\zo/{$0\mkern2mu\mhyphen1$}

\def\nn/{$n \times n$}

\title{Random Regular Bipartite Graphs Satisfy Weak Virial Positivity, for a Large Range of the Parameters}
\author{Paul Federbush\\
Department of Mathematics\\
University of Michigan\\
Ann Arbor, MI, 48109-1043\\
pfed@umich.edu}

\date{\today}

\usepackage{amsthm} 
\newtheorem{thm}{Theorem}[section]
\newtheorem{conj*}{Conjecture} 

\numberwithin{lemma}{section}
\numberwithin{conj*}{section}
\DeclareMathOperator{\Prob}{\mathrm{Prob}}
\DeclareMathOperator{\EE}{\mathbf{E}}

\numberwithin{equation}{section}

\begin{document}

\maketitle

\section{Introduction and Results}\label{sec:I}
We deal with $r$-regular bipartite graphs with $v = 2n$ vertices. Let $m_i$ be the number of $i$-matchings. In \cite{1}, Butera, Pernici, and I introduced the quantity $u(i)$ in eq (1) therein, 
\begin{align}\label{I.1}
	u(i) = -\ln(i!m_i)
\end{align}
but we did not use this notation there. We then considered $\Delta^k(u(i))$ where $\Delta$ is the finite difference operator so
\begin{align}\label{I.2}
	\Delta u(i) = u(i+1) - u(i).
\end{align}

A graph is said to satisfy Virial Positivity if all the meaningful 
$\Delta^k u(i)$, $k \geq 2$, were non-negative. That is
\begin{align}\label{I.3}
	\Delta^k u(i) \geq 0
\end{align}
for $k = 2,\dots, v$ and $i = 0,\dots, v-k$. We make the ``Virial Positivity Conjecture" supported by some computer evidence.

\begin{conj*}\label{conj*:I.1}
As $n$ goes to infinity the fraction of graphs that satisfy Virial Positivity approaches $1$.
\end{conj*}

We note some of the impressive results of the numerical study of Virial Positivity in \cite{1}. 
\begin{enumerate}
	\item All graphs with $v < 18$ satisfy Virial Positivity.
	\item When $r = 4$, the first violation occurs when $v = 20$ in a single graph among 62611 graphs with $v = 20$.
	\item For $r = 3$ the fraction of graphs not satisfying Virial Positivity monotonically decreases between 18 and 30 with a single violation at 18. 
\end{enumerate}

The term Virial comes from the Virial expansion in statistical mechanics. Equation (\ref{I.1}) above corresponds, roughly speaking,  to the positivity of the coefficients in the Virial Expansion for infinite regular lattices, \cite{1}. 
For hyper-cubical lattices it is shown in \cite{1} that the first 20 coefficients in the Virial expansion are positive in dimensions 
$d \leq 10$ !

Our treatment of Virial Positivity follows slavishly the treatment of Graph Positivity in \cite{3}. Graph Positivity was first defined in \cite{2}, generalizing the positivity of certain coefficients on infinite regular lattices computed in \cite{6}. 

The results in the present paper are basically corollaries of parallel results in \cite{3} for Graph Positivity. To understand the results one need not consult \cite{3}, but when we discuss the proof, in Section \ref{sec:II}, knowledge of the treatment in \cite{3} is indispensable. And we do not claim that \cite{3} is not a difficult paper. As the double negative implies, the developments in \cite{3} require serious study for grasp.

In this paper we study a weaker property than the Virial Positivity conjecture. We work with the Weak Virial Positivity Conjecture.

\begin{conj*}\label{conj*:I.2}
For each $i$ and $k \geq 2$ one has 
\begin{equation*}
\mathrm{Prob}(\Delta^k u(i) \geq 0) \xrightarrow[n \to \infty]{} 1.
\end{equation*}
In fact, what we prove is
\end{conj*}

\begin{thm}\label{thm:I.1}
If $r \leq 10$, $i + k \leq 100$, $2 \leq k \leq 27$, or $i + k \leq 29$ all r, then
\begin{equation*}
\mathrm{Prob}(\Delta^k u(i) \geq 0) \xrightarrow[n \to \infty]{} 1.\end{equation*}
\end{thm}

The proof is presented in Section \ref{sec:II}. We here note that one relies heavily on the work of Wanless, \cite{4}, and Pernici, \cite{5}. We are also indebted to Robin Chapman, see Appendix C of \cite{3}, for a combinatoric proof. The extension of the range of validity of the
theorem to  $i + k \leq 29$ all $r$ is due to numerical computer
computations presented in \cite{7}.

There is much room for future work: extending the range of parameters in this paper and \cite{3}, finding the status of the strong forms of the conjectures, understanding of the positivities alluded to on infinite regular lattices. Above all, why these positivities!

\section{Proof}\label{sec:II}
The proof of Theorem \ref{thm:I.1} is constructed by assembling pieces extracted from the developments in \cite{3}, making small modifications. No additional real cleverness or hard computations are required. Strangely enough most of our effort is in adapting the material of Section 2, the easiest section in \cite{3}.

\subsection*{An Avatar of Section 2 of \cite{3}}
Suppose we want $\mathrm{Prob}(x<y)$ to be large. We have 
\begin{align}\label{II.1}
	\Prob(x > y) = \Prob(e^x > e^y).
\end{align}
Set 
\begin{align}
	e^x - e^y \equiv \alpha_0
\end{align}
and 
\begin{align}
	\EE(e^x - e^y) \equiv \alpha.
\end{align}
We will want $\alpha$ to be negative.
\begin{align}
	\EE((e^x-e^y)^2) - \alpha^2 \equiv \beta.
\end{align}
Then assuming $\alpha$ is negative,
\begin{align}
	\EE((e^x - e^y - \alpha)^2) \geq \alpha^2 \Prob(e^x -e^y > 0).
\end{align}
And so
\begin{align}\label{II.6}
	\Prob(e^x > e^y) \leq \dfrac{\beta}{\alpha^2}.
\end{align}

In our problem $\beta$ and $\alpha$ will be functions of $n$ and we'll want probability to go to zero with $n$ as $n$ goes to infinity. 

We turn to the object of study of eq. (\ref{I.1}).
\begin{subequations}
\begin{align}\label{II.7a}
	\Prob(\Delta^k u(i) < 0) = \Prob\left((-1)^k \sum_{l=0}^k (-1)^k \binom{k}{l} \ln\left(m_{i+l}(i+l)!\right) <0\right).
\end{align}

A deus ex machina now imports eq. (3.11) of \cite{3} the relation
\begin{align*}
	m_i = \dfrac{n^ir^i}{i!} (1+ \hat {H_i })
\end{align*}
and substitutes it into eq. (\ref{II.7a}) 
\begin{align}
	= \Prob\left((-1)^k\sum_{l=0}^k (-1)^k \binom{k}{l} \left[ \ln(1+ \hat {H_{i+l}}) + (i+l) \ln(nr)\right]<0\right)
\end{align}
which becomes if $k \geq 2$
\begin{align}
	= \Prob\left((-1)^k \sum_{l=0}^k (-1)^k \binom{k}{l} \ln\left(1+\hat {H_{i+l}} )\right) <0\right) 
\end{align}
\end{subequations}
\begin{align}
	 = \Prob\left(\sum_{l \in \mathcal{L}^+} \binom{k}{l} \ln(1+\hat {H_{i+l}} ) < \sum_{l \in \mathcal{L}^-} \binom{k}{l} \ln(1+\hat {H_{i+l}} )\right)
\end{align}
where $\mathcal{L}^+$ is the set of odd $l$,$ 0 \leq l \leq k$, if $k$ is odd, and is the set of even $l$, $0 \leq l \leq k$ if $k$ is even and $\mathcal{L}^-$ is defined vice versa.

Referring to (\ref{II.1}) through (\ref{II.6}) we set

\begin{align}
	x &= \sum_{l \in \mathcal{L}^+} \binom{k}{l} \ln (1+ \hat {H_{i+l}})\\
	y &= \sum_{l \in \mathcal{L}^-} \binom{k}{l} \ln (1+ \hat {H_{i+l}})\\
	e^x &= \prod_{l \in \mathcal{L}^+} (1+ \hat {H_{i+l}})^{\binom{k}{l}
	}\\
	e^y &= \prod_{l \in \mathcal{L}^-} (1+ \hat {H_{i+l}})^{\binom{k}{l}}
\end{align}
and
\begin{align}\label{II.13}
	\alpha_0 = \prod_{l \in \mathcal{L}^+} (1+ \hat {H_{i+l}})^{\binom{k}{l}} - \prod_{l \in \mathcal{L}^-} (1+ \hat {H_{i+l}})^{\binom{k}{l}}.
\end{align}

\subsection*{Extracting the Proof From \cite{3}}
The proof we want of Theorem I.1 can be obtained from the results of \cite{3} by the following wonderfully simple procedure:

We compare eq. (\ref{II.13}) above with eq. (4.2) of \cite{3} which we rewrite here

\begin{equation}\label{eq4.2}\tag*{(4.2)}
 \alpha_0=\Biggl(
    \prod_{\ell \in \mathcal{L}^+} \bigl(
      (1+\hat{H_{i+\ell}})(1+K_{i+\ell})
    \bigr)^{\binom{k}{\ell}}
    - \prod_{\ell \in \mathcal{L}^-} \bigl(
      (1+\hat{H_{i+\ell}})(1+K_{i+\ell})
    \bigr)^{\binom{k}{\ell}}
  \Biggr)
\end{equation}

SO $\alpha_0$ OF THIS PAPER IS THE SAME AS THE  $\alpha_0$ of \cite{3} WITH THE $K_i$ ALL SET EQUAL TO  $0$  !  

Whereas \cite{3} is a difficult paper it is easy to go through and set $K_i = 0$ in the essential places. We proceed to do so. From the discussion between eq. (6.9) and eq. (6.11) we change the right side of (6.11) 
\begin{align}\label{new6.11}\tag*{new (6.11)}
	= \dfrac{(k-2)!}{k!}\left(\dfrac{1}{r^{k-1}} -2\right)i^k
\end{align}
We change then (6.5)
\begin{align}\label{new6.5}\tag*{new (6.5)}
	\dfrac{(k-2)!}{k!}\left(\dfrac{1}{r^{k-1}}-2\right)i^k
\end{align}
and then the right side of (6.3)
\begin{align}\label{new6.3}\tag*{new (6.3)}
	=(k-2)!\left(\dfrac{1}{r^{k-1}}-2\right)
\end{align}
and following, right side of (8.7)
\begin{align}\label{new8.7}\tag*{new (8.7)}
	=\begin{cases} 0 & \qquad d < k-1 \\
		(k-2)!\left(\dfrac{1}{r^{k-1}}-2\right) &\qquad d=k-1
	\end{cases}
\end{align}
and finally (8.1)
\begin{align}\tag*{new (8.1)}
\alpha_0 = \dfrac{(k-2)!}{n^{k-1}}\left(\dfrac{1}{r^{k-1}}-2\right) + \mathcal{O}\left(\frac{1}{n^k}\right)
\end{align}

Noting that $\alpha_0$ is negative and using the argument at the end of Section 9 of \cite{3} we have a proof of Theorem \ref{thm:I.1}. We hope the reader is inclined to study \cite{3}, there is a lot there. 

As a final note, it is possible to prove by the same method the conjecture made in \cite{1} that

$$
\mathrm{Prob}(\Delta^k (-m_i i!(n-i)!) \geq 0) \xrightarrow[n \to \infty]{} 1, \quad k\geq 2.
$$

The thoughtful reader may come to believe that there are more such results provable than one would desire.

\end{document}